\documentclass[draft]{amsart}
\usepackage{amsmath,amsthm,amssymb}
\newtheorem{thm}{Theorem}
\newtheorem{lem}{Lemma}
\newtheorem{cor}{Corollary}

\theoremstyle{remark}
\newtheorem{example}{Example}
\newcommand{\C}{\mathbb{C}\mkern1mu}
\newcommand{\N}{\mathbb{N}\mkern1mu}

\newcommand{\PP}{\mathbb{P}\mkern1mu}
\newcommand{\cond}{\,\vert\,}

\newcommand{\bmat}{\begin{pmatrix}}
\newcommand{\emat}{\end{pmatrix}}
\newcommand{\bsmat}{\bigl(\begin{smallmatrix}}
\newcommand{\esmat}{\end{smallmatrix}\bigr)}
\newcommand{\Bsmat}{\Bigl(\begin{smallmatrix}}
\newcommand{\Esmat}{\end{smallmatrix}\Bigr)}
\newcommand{\bbsmat}{\biggl(\begin{smallmatrix}}
\newcommand{\eesmat}{\end{smallmatrix}\biggr)}
\begin{document}
\title{Chow points of $\C$-orbits}
\author{Annett P\"uttmann}
\date{}
\thanks{Research supported by SFB/TR12 ``Symmetrien und Universalit\"at in mesoskopischen 
Systemen'' of the Deutsche Forschungsgemeinschaft.}
\address{Fakult\"at f\"ur Mathematik, Ruhr-Universit\"at Bochum, 44780 Bochum}
\email{annett.puettmann@rub.de}
\begin{abstract}
We consider free algebraic actions of the additive group of complex numbers on a 
complex vector space $X$ embedded in the complex projective space.
We find an explicit formula for the map $p$
that assigns to a generic point $x\in X$ the Chow point of the closure of the orbit
through $x$.
The properties Hausdorff quotient topology and proper action are equivalently 
characterized by the closure of the image of $p$ in the closed Chow variety.
\end{abstract}
\maketitle
\section{introduction}
In this article we consider free algebraic $\C$-actions, i.e., 
actions of the additive group of complex numbers, on $\C^n=X$.
Such an action is given by a locally nilpotent derivation $\delta:\C[X|\to\C[X]$.
We have $t.f=\sum_{j=0}^{\infty} \frac{t^j}{j!}\delta^j f$ for $t\in\C$ and $f\in\C[X]$, where
$\C\times\C[X]\to\C[X]$, $(t,f)\mapsto t.f$, denotes the action of $\C$ on $\C[X]$.

Closures of $\C$-orbits in a compactification $\PP_n\C$ of $X$ are projective curves in
$\PP_n\C$.
In section \ref{sec:ChowPoint}, Theorem \ref{thm:ChowPointGenericOrbit}, 
we derive an explicit formula for the map 
$p: U_0 \to {\mathcal C}_{d,n}$ from the Zariski-open set $U_0 \subset X$ of orbits of 
maximal degree $d$ to the open Chow variety ${\mathcal C}_{d,n}$ of curves of degree $d$
in $\PP_n\C$ that assigns to a point $x\in X$ the Chow point of 
$\overline{\C.x} \subset \PP_n\C$.

Since the closed Chow variety $\overline{{\mathcal C}_{d,n}}$ is projective,
there is a sequence of blow-ups $\pi:\tilde{X} \to X$ and an regular map 
$\tilde{p}: \tilde{X} \to \overline{{\mathcal C}_{d,n}}$ that lifts $p$.

We are interested in properties of the $\C$-action on $\C^n$ that
are relevant for geometric quotients.
By definition, the quotient topology is Hausdorff if the orbits can be separated by
continuous functions. 
Recall that the action of a Lie group $G$ on a manifold $X$ is called proper if the map
$G\times X\to X\times X$, $(g,x)\mapsto(x,g.x)$, is proper.

In section \ref{sec:LimitCycles} we express  {\it properness}  and 
{\it Hausdorff quotient topology} using properties of the cycles in $\tilde{p}(\tilde{X})$.
Corollary \ref{cor:hausdorff} states that the quotient topology is Hausdorff if and only if 
for all $\tilde{x} \in \tilde{X}$ 
all irreducible components of $\tilde{p}(\tilde{x})$ except the closure of $\C.\pi(\tilde{x})$
are contained in infinity.
In Corollary \ref{cor:proper} we prove that the action is proper if and only if it has Hausdorff
quotient topology and the multiplicity of the closure of $\C.\pi(\tilde{x})$ in
$\tilde{p}(\tilde{x})$ is $1$ for all $\tilde{x} \in \tilde{X}$.

We exploit the fact that the analytic space associated the Chow scheme of curves in $\PP_n\C$
is the Barlet space of compact $1$-cycles in $\PP_n\C$ \cite{B}.
This enables us to describe 
convergence of Chow points by the behavior of the corresponding cycles.

In section \ref{sec:examples} we apply our results to free $\C$-actions of degree two
that have been studied by other authors using different methods.
So far our approach does not describe the geometric quotient in the case 
where the action is proper.
But the map $\tilde{p}$ gives complete information about non-separable orbits if the
quotient topology is not Hausdorff.

The construction of the explicit form of the map $p: U_0 \to {\mathcal C}_{d,n}$ strongly
depends on the fact that the acting group is $\C$.
Whereas it should be possible to generalize the results of section \ref{sec:LimitCycles} 
to free actions of other groups on more general affine varieties than $\C^n$.

We would like to thank A.~T.~Huckleberry for valuable discussions.
\section{Chow point of a generic orbit}
\label{sec:ChowPoint}
Fixing a set of coordinates $\{ x_1,\ldots,x_n\}\subset\C[X]$ on $X=\C^n$ defines an embedding
$\C^n\hookrightarrow\PP_n\C$ by $x\mapsto [1:x_1(x):\ldots:x_n(x)]$.
Furthermore, every derivation $\delta$ can be uniquely written as 
$\delta=\sum_{j=1}^n \delta(x_j)\frac{\partial}{\partial x_j}$.
The $\C$-action on $X$ corresponding to a locally nilpotent derivation $\delta$ on $\C[X]$
is  completely described by $x_k(-t.x)=\sum_{j=0}^d \frac{t^j}{j!}(\delta^jx_k)(x)$.

For a locally nilpotent derivation $\delta$ on $\C[X]$ we define its {\it degree} to be the
smallest number $d$ such that $\delta^{d+1}(x_j)=0$ for all $j$.
Of course, this notion of degree depends on the chosen coordinates.
Let  $C_x:=\overline{\C.x}\subset\PP_n\C$ be the closure in $\PP_n\C$ of the orbit through 
the point $x\in X$.
Although all orbits of a free $\C$-action are isomorphic to $\C$ their closures in $\PP_n\C$
may lead to non-isomorphic projective curves distinguishable  by their degree as 
projective varieties.
\begin{lem}
If $\delta$ defines a free algebraic $\C$-action on $\C^n$ of degree $d$, 
then $\deg C_x\leq d$  for all $x\in\C^n$.
Furthermore, $\deg C_x=d$ iff $x\in\bigcup_{j=1}^n \{ \delta^d(x_j)\ne 0\}$.
\end{lem}
\begin{proof}
For a general hyperplane $H=\{ \sum_{j=0}^n a_jX_j=0\}$ in $\PP_n\C$ the intersection
$H\cap C_x$ is given by the solutions of
$a_0+\sum_{k=1}^n a_k \sum_{j=0}^d \frac{t^j}{j!}(\delta^jx_k)(x) =0$
which is a polynomial in $t$ of degree at most $d$.
\end{proof}
The set $U_0:=\cup_{j=1}^n \{ \delta^d(x_j)\ne 0\}$ is Zariski-open and $\C$-invariant, because
$\delta^{d+1}x_j=0$ for all $j$ and $\delta^d x_j\ne 0$ for at least one $j$.
An orbit in $U_0$ is called {\it generic orbit}.
Hence, given a free algebraic $\C$-action on $\C^n$ of degree $d$ and a compactification
$\C^n \subset \PP_n\C$
the closures of generic orbits are projective curves of degree $d$ in $\PP_n\C$.
Projective curves of degree $d$ in $\PP_n\C$ are parametrized by the open
Chow variety ${\mathcal C}_{d,n}$. 
We will give an explicit description of the map $p:U_0\to{\mathcal C}_{d,n}$ that
assigns to a point $x\in X$ the Chow point of $C_x$.

We refer to \cite{H} for an introduction to the Chow variety.
The Chow variety ${\mathcal C}_{d,n}$ is a quasiprojective variety in 
$\PP(\Gamma(\PP_n^2(\C),H^d))$.
Here, $\Gamma(\PP_n^2(\C),H^d)$ is the vector space of bihomogeneous polynomials 
of bidegree $(d,d)$  with respect to the two sets of $n+1$ variables $\alpha_0,\ldots,\alpha_n$ 
and $\beta_0,\ldots,\beta_n$.
The {\it Chow point} $P_C\in{\mathcal C}_{d,n}$ of a curve $C$ of degree $d$ in 
$\PP_n\C$ is the polynomial defining the hypersurface  
$H_C:=\cup_{x\in C}\{ (\alpha,\beta)\cond x\in H_{\alpha}\cap H_{\beta}\}\subset\PP_n^2(\C)$.
Here, the hyperplane $H_{\alpha}\subset\PP_n\C$ assigned to a point 
$\alpha\in\PP_n\C$ is $H_{\alpha}=\{\sum_{j=0}^n \alpha_jX_j = 0\}$.

We denote by $\tau:\PP_n^2(\C)\to\PP_n^2(\C)$, $(\alpha,\beta)\mapsto(\beta,\alpha)$, the
permutation of the components of $\PP_n^2(\C)$.
Note that $\gamma_{j_1j_2}:=\beta_{j_1}\alpha_{j_2}-\beta_{j_2}\alpha_{j_1}$
are $\tau$-invariant as elements of $\PP(\Gamma(\PP_n^2(\C),H^1))$ for $j_1,j_2=0,\ldots,n$.
For $k,l=0,\ldots,d$ we set
\begin{equation}
\label{eq:flk}
 f_{lk}(x) := 
 \sum_{j=0}^n\beta_j(\delta^lx_j)(x)\sum_{j=0}^n\alpha_j(\delta^kx_j)(x)
-\sum_{j=0}^n\alpha_j(\delta^lx_j)(x)\sum_{j=0}^n\beta_j(\delta^kx_j)(x)
\end{equation}
Note that $f_{kl}(x) = \sum_{j_1=0}^n\sum_{j_2=0}^n \gamma_{j_1j_2}
\delta^l(x_{j_1})(x)\delta^k(x_{j_2})(x)$.
Since  $\gamma_{j_1j_2} = \tau.\gamma_{j_1j_2} = -\gamma_{j_2j_1}$, we have
$f_{lk}=-f_{kl}$. Let
\begin{align}
\label{eq:Flk}
F_{lk}(x)  & := \sum_{r=0}^{\min{d-l,k}} f_{l+r, k-r}(x)
	\prod_{j=0}^{d-l-r-1}(d-j)\prod_{j=0}^{r-1} (k-j), \\
\label{eq:P(x)}
P(x) & := \det( F_{lk}(x) )_{l\ne 0,k\ne d}
\end{align}
\begin{lem}
The polynomial $P$ is a $\tau$-invariant element of $\PP(\Gamma(\PP_n^2(\C),H^d))$.
\end{lem}
\begin{proof}
By construction, $P$ is a homogeneous polynomial of degree $d$ in the terms
$\gamma_{lk}$.
It remains to show that $P$ is non-trivial, $P\not\equiv 0$.

To see this, we examine the coefficients $c_{j0}$ of the monoms $\gamma_{j0}^d$.
Among the polynomials $f_{lk}$ only $f_{0k}$ and $f_{l0}$
contribute to the coefficients $c_{j0}$, because $\delta^l(x_0)=0$ for all $l>0$.
Thus, the coefficient of $\gamma_{j0}$ in $F_{lk}$ vanishes if $d-l < k$.
Since the definition of $P$ involves a determinant, $c_{j0}$ equals the 
coefficient of $\gamma_{j0}^d$ in $f_{d0}^d$.
By definition, $f_{d0}= (\sum_{j=1}^n\gamma_{j0}\delta^d x_j+
\sum_{j_1=1}^n\sum_{j_2=1}^n \gamma_{j_1j_2} (\delta^d x_{j_1})x_{j_2})(x)$.
Hence, $c_{j0}=(\delta^d x_j)^d(x)$. 
Since the action has degree $d$, we have $c_{j0}\ne 0$ for at least one $j$.
\end{proof}
\begin{lem}
\label{lem:fdFd-1}
Let $f_{lk}$ be elements of a polynomial ring
satisfying $f_{lk}=-f_{kl}$ for all $k,l=0,\ldots,d$ and $t\in\C$.
The equations $\sum_{k=0}^d f_{lk}\frac{t^k}{k!}=0$ are satisfied 
for all $l=0,\ldots,d$ iff
$\sum_{k=0}^{d-1} F_{lk} \frac{t^k}{k!}=0$ for all $l=1,\ldots,d$, where
$F_{lk}$ is defined as in (\ref{eq:Flk}).
\end{lem}
\begin{proof}
Consider the following recursion:
$f^{(0)}_{lk}:=f_{lk}$, 
\[ f^{(s+1)}_{lk} := \begin{cases} f^{(s)}_{lk} & l\ne d-s-1 \\
	f^{(s)}_{d-s, k-1}k+f_{d-s-1, k}\prod_{j=0}^s(d-j) & l=d-s-1\end{cases}\]
We will prove by induction on $s$ that $f^{(s)}_{ld}=0$ for $l\geq d-s$ and
\[ f^{(s)}_{lk} = \begin{cases} f_{lk} & l<d-s \\
	\sum_{r=0}^{\min{d-l,k}} f_{l+r, k-r}
 \prod_{j=0}^{d-l-r-1}(d-j)\prod_{j=0}^{r-1}(k-j) & 
		l\geq d-s\end{cases} \]

If $s=0$, then $f_{dk}^{(0)}= f_{dk} = \sum_{r=0}^0 f_{d+0, k-0}$. In particular,
$f_{dd}^{(0)}=0$. Now,
\begin{align*}
f^{(s+1)}_{d-s-1, k} =& k\sum_{r=0}^{\min{s,k-1}} f_{d-s+r, k-1-r} 
	\prod_{j=0}^{s-r-1}(d-j)\prod_{j=1}^r(k-j)+f_{d-s-1, k}
	 \prod_{j=0}^s(d-j)\\
	=& \sum_{r=1}^{\min{s+1,k}} f_{d-s-1+r, k-r}
	 \prod_{j=0}^{s-r}(d-j)\prod_{j=0}^{r-1}(k-j)
		+f_{d-s-1, k}\prod_{j=0}^s(d-j)\\
	= & \sum_{r=0}^{\min{s+1,k}} f_{d-s-1+r, k-r}
	 \prod_{j=0}^{s+1-r-1}(d-j)\prod_{j=0}^{r-1}(k-j)
\end{align*}
Furthermore, $f^{(s+1)}_{d-s-1,d}=\sum_{r=0}^{s+1} f_{d-s-1+r, d-r}
\prod_{j=0}^{s+1-r-1}(d-j)\prod_{j=0}^{r-1}(d-j) =0$, 
since $f_{d-(s+1-r), d-r}=-f_{d-r, d-(s+1-r)}$.

Multiplying the equation $\sum_{k=0}^d f^{(s)}_{d-s, k}\frac{t^k}{k!}=0$ by $t\ne 0$
gives an equivalent equation
$\sum_{k=0}^{d-1} f^{(s)}_{d-s, k}\frac{t^{k+1}}{k!}= 
\sum_{k=0}^d f^{(s)}_{d-s, k-1}k\frac{t^k}{k!}=0$, since $f^{(s)}_{d-s, d}=0$.
Consequently, the systems  $\sum_{k=0}^d \frac{t^k}{k!} f^{(s)}_{l k}=0$ for
all $l=0,\ldots,d$ are equivalent for all $s$. 

Note that $F_{lk}=f^{(d)}_{lk}$. 
Recall that $F_{ld}=0$ for all $l$.
Since $f_{r, k-r}=-f_{k-r, r}$, we have for all $k$
\[ F_{0k}=\sum_{r=0}^k f_{r, k-r}\prod_{j=0}^{d-r-1}(d-j)\prod_{j=0}^{r-1} (k-j) 
	= \sum_{r=0}^k f_{r, k-r}\tfrac{d!}{r!}\tfrac{k!}{(k-r)!}=0.\]
\end{proof}
\begin{thm}
\label{thm:ChowPointGenericOrbit}
Let $\delta$ be a locally nilpotent derivation of $\C[x_1,\ldots,x_n]$
that defines a free $\C$-action on $\C^n$ of degree $d$ and $x\in \C^n$.
If $\C.x$ a generic orbit, then the Chow point $p(x)$ of $C_x$ is 
$P(x)$ as defined by equation (\ref{eq:P(x)}).
\end{thm}
\begin{proof}
Except for one point, the point at infinity, the curve $C_x$ is parametrized by the
$\C$-action.
Therefore, $H_{C_x}=\overline{\{ (\alpha,\beta)\cond 
	\exists t\in\C: t.x\in H_{\alpha}\cap H_{\beta} \}}$.
To find the defining polynomial of of $H_{C_x}$, we can assume that
if $(\alpha,\beta)\in H_{C_x}$ then there is $t\in\C$ such that
$\sum_{j=0}^n \alpha_j x_j(-t.x)=0$ and $\sum_{j=0}^n \beta_j x_j(-t.x)=0$.
Since for all $j$ $t.x_j=\sum_{k=0}^d \frac{t^k}{k!}\delta^k x_j$, 
we obtain the equations
\begin{eqnarray}
\label{eq:Ak}
\sum_{k=0}^d \tfrac{t^k}{k!}A_k=0, \quad A_k:=\sum_{j=0}^n\alpha_j\delta^k x_j(x)\\ 
\label{eq:Bk}
\sum_{k=0}^d \tfrac{t^k}{k!}B_k=0, \quad B_k:=\sum_{j=0}^n\beta_j\delta^k x_j(x)
\end{eqnarray}
The terms $A_k$ and $B_k$ are polynomials of bidegree $(1,0)$ and $(0,1)$,
respectively.
The combination of the equations (\ref{eq:Ak}) and (\ref{eq:Bk}) leads to a system of $d+1$ equations 
\begin{equation}
\label{eq:AB}
\sum_{k=0}^d\frac{t^k}{k!}(B_lA_k-A_lB_k)=0,\qquad l=0,\ldots,d.
\end{equation}
Note that $f_{lk}=B_lA_k-A_lB_k$.
By Lemma \ref{lem:fdFd-1}, the system (\ref{eq:AB}) is equivalent to
\begin{equation*}
\sum_{k=0}^{d-1}\frac{t^k}{k!}F_{lk} = 0, \qquad l=1,\ldots,d
\end{equation*}
Since $t^0=1$, it follows that $P=\det(F_{lk})_{l\ne 0,k\ne d}=0$
for all $(\alpha,\beta)\in H_{C_x}$.
But $P$ is a homogeneous polynomial of bidegree $(d,d)$ and $\tau$-invariant
as an element of $\PP(\Gamma(\PP_n^2(\C),H^d))$.
\end{proof}
A {\it local slice} of an action of a Lie group $G$ on a manifold $X$ 
through a point $x\in X$ is a locally closed  submanifold $x\in S\subset X$ such that 
$G.S$ is open in $X$ and  every orbit through $S$ intersects $S$ in exactly one point.
For a free action a local slice establishes a local trivialization of the action, i.e.,
$G\times S \to G.S$, $(g,s) \mapsto g.s$, is a diffeomorphism.
The properness of a Lie group action implies the existence of local slices.
We include the proof of a special case of this general fact.
\begin{lem}
A free action of a Lie group $G$ on a manifold $M$ is proper if and only if
there are local slices through every point of $M$.
A  proper action of a Lie group $G$ on a manifold $M$has Hausdorff quotient topology.
\end{lem}
\begin{proof}
The action of a Lie group $G$ on a manifold $M$ is proper iff 
for any convergent sequence $\{ x_n \} \subset M$  and any sequence $\{ g_n\} \subset G$
the convergence of the sequence $\{ g_n.x_n\} \subset M$ implies the existence of a convergent 
subsequence of $\{ g_n\} \subset G$.

If there is a local slice through every point in $M$, then $M$ can be covered by local 
trivializations of the action. Then a local trivialization around $y:= \lim g_n.x_n$ shows that
$\{ g_n\}$ is convergent.

Assume that $G$ acts properly on $M$. Let $x\in M$ be a point.
Choose any locally closed submanifold $S'$ through $x$ 
that is transversal to the $G$-orbit $G.x$.
Such a $S'$ can be shrunk to $S$ such that $G\times S\to G.S$, $(g,s)\mapsto g.s$, is 
a diffeomorphism. Otherwise two convergent sequences $\{ s_n\}\to x$, $\{s'_n=g_n.s_n\}\to x$
would exist. The properness would imply that $\{ g_n\}$ is convergent to the identity $e$
contradicting the fact that $G\times S\to G.S$ is locally diffeomorphic at $(e,x)$. 

If $\{ x_n\}$ and $\{ g_n.x_n\}$ are convergent sequences in $M$, then the properness of
the action implies that $g_n$ tends to an element $g\in G$. 
Then $g.\lim_{n\to \infty} x_n = \lim_{n\to \infty} g_n.x_n$. 
\end{proof}
\begin{lem}
There are local slices through every $x\in U_0$.
\end{lem}
\begin{proof}
There are local trivializations 
$\{ \delta^d x_j \ne 0\}\to \C\times\{\delta^d x_j\ne 0, \delta^{d-1}x_j =0  \}$ given by
$x\mapsto (t(x), -t(x).x)$ with $ t(x) = (\delta^{d-1}x_j)(x) / (\delta^dx_j)(x) $,
because 
\[Ê \delta^{d-1}x_j(-t(x).x)=\sum_{k=0}^d \frac{(-t(x))^k}{k!}\delta^k\delta^{d-1}x_j(x)
=\delta^{d-1}x_j(x)-t(x)\delta^d x_j(x) = 0. \]
\end{proof}
\begin{cor}
The restriction of the $\C$-action on $U_0$ is a proper action.
\end{cor}
Note that $U_0 = X$ if $d=1$. Hence, a free algebraic $\C$-action of degree one 
on $\C^n$ is proper.
\section{Limit cycles }
\label{sec:LimitCycles}
The regular map $p:U_0\to{\mathcal C}_{d,n}$, $x\mapsto P(x)$ defines a rational map  
$p:X\to{\mathcal C}_{d,n}$.
Since the open Chow variety ${\mathcal C}_{d,n}$ is quasiprojective, 
there is a finite sequence of blow-ups
$\pi:\tilde{X}=X_m\to X_{m-1}\to\ldots\to X_0=X$ along $Y_j\subset X_j$ 
and a regular morphism $\tilde{p}:\tilde{X}\to\overline{{\mathcal C}_{d,n}}$
into the projective Chow variety $\overline{{\mathcal C}_{d,n}}$ 
such that $\tilde{p}=p\circ\pi$.
A point in $\overline{{\mathcal C}_{d,n}}$ is a cycle $C=\sum_{j=1}^s n_jC_j$ which consists
of irreducible curves $C_j \subset\PP_n \C$ with multiplicities $n_j \in \N$ such that
$d=\sum_{j=1}^s n_j\deg C_j$.

For $\tilde{x} \in \tilde{X} \setminus \pi^{-1}(U_0)$ we call 
$\tilde{p}(\tilde{x}) \in {\mathcal C}_{d,n}$ a {\it limit cycle}.
A cycle  $\tilde{p}(\tilde{x})$ can be uniquely written as
$\tilde{p}(\tilde{x})=n_{\tilde{x}}C_{\tilde{x}}+Z_{\tilde{x}}$ with $n_{\tilde{x}}\in\N$,
$C_{\tilde{x}}=\overline{\C.\pi(\tilde{x})}$, and a cycle $Z_{\tilde{x}} \subset \PP_n\C$
of degree $d-n_{\tilde{x}} \deg C_{\tilde{x}}$ that does not contain $C_{\tilde{x}}$.
Note that $n_{\tilde{x} }= 1$ and $Z_{\tilde{x}} = 0$ if $\pi(\tilde{x}) \in U_0$.

In the sequel we us the fact that the analytic space associated to the closed Chow variety
is contained in the Barlet space.
\begin{thm}
If $Y$ is a projective variety, then the Barlet space of $Y$ coincides with the analytic
space associated with the Chow scheme of $Y$.
\end{thm}
\begin{proof}
See \cite{B} and \cite{M}.
\end{proof}
This allows us to discuss convergence of Chow points $\tilde{p}(x_n)$ in terms of the 
corresponding projective curves $C_{x_n}$.
\begin{lem}
\label{lem:pointwise}
Let $\tilde{x} \in \tilde{X}$.
For every open neighborhood $W\subset \PP_n\C$ of $\tilde{p}(\tilde{x})$, there is an  open
neighborhood $U\subset \tilde{X}$ of $\tilde{x}$ such that
$\tilde{p}(y) \subset W$ for all $y\in U$.
\end{lem}
\begin{proof}
The family of compact $1$-cycles parametrized by the points in the Barlet space 
is an analytic family of compact $1$-cycles, see \cite{B} and \cite{M}.
\end{proof}
\begin{cor}
For all $\tilde{x}\in\tilde{X}$ the curve $C_{\tilde{x}}$ is a component of the limit cycle 
$\tilde{p}(\tilde{x})$, i.e., $n_{\tilde{x}}\geq1$.
\end{cor}
\begin{proof}
Let $\tilde{x}\in\tilde{X}$. We take a sequence $\{ \tilde{x}_n\}\subset \pi^{-1}(U_0)$
that tends to $\tilde{x}$ in $\tilde{X}$.
Now, $\tilde{p}(\tilde{x}_n) \to \tilde{p}(\tilde{x})$, since $\tilde{p}$ is continuous.
But $\tilde{p}(\tilde{x}_n) = p(\pi(\tilde{x}_n))$ and the sequence $\{\pi(\tilde{x}_n)\}$ 
tends to $\pi(\tilde{x})$.
By lemma \ref{lem:pointwise}, $C_{\tilde{x}}$ is a component of $\tilde{p}(\tilde{x})$, 
because $t.\pi(\tilde{x}_n) \to t.\pi(\tilde{x})$ for all $t \in \C$.
\end{proof}
\begin{cor}
\label{cor:hausdorff}
The $\C$-action on $X$ has Hausdorff quotient topology if and only if
$Z_{\tilde{x}} \cap X = \emptyset$ for all $\tilde{x}\in\tilde{X}$.
\end{cor}
\begin{proof}
Two distinct orbits $\C.x$ and $\C.x'$ in $X$  can not be separated by  
$\C$-invariant open subsets if and only if for any pair of open neighborhoods $U$ of $x$ and
$U'$ of $x'$ there exists a $t\in \C$ such that $t.U \cap U' \ne \emptyset$. 
This means that there are convergent sequences
$\{\tilde{x}_n\}\to\tilde{x}$, $\{\tilde{x}'_n\}\to\tilde{x}'$ in $\tilde{X}$ and a sequence
$\{ t_n \} \subset \C$ such that $\pi(\tilde{x}'_n) = t_n.\pi(\tilde{x}_n)$,
$\pi(\tilde{x}_n)$ tends to $x$, and $\pi(\tilde{x}'_n)$ tends to $x'$.
By lemma \ref{lem:pointwise} both limit cycles, $\tilde{p}(\tilde{x})$ and $\tilde{p}(\tilde{x}')$,
contain the irreducible components $C_{\tilde{x}}$ and $C_{\tilde{x}'}$.
\end{proof}
The condition $Z_{\tilde{x}} \cap X = \emptyset$ means that $Z_{\tilde{x}}$
contains only points at infinity, i.e., $Z_{\tilde{x}} $ is contained in the hyperplane 
$\{ X_0 = 0 \} \subset \PP_n\C$.
Hence,  $Z_{\tilde{x}} \cap X = \emptyset$ iff 
$\tilde{p}(\tilde{x})=p(\pi(\tilde{x}))^{n_{\tilde{x}}}r$ where $r$ is a polynomial in the set of 
variables $\{Ê\gamma_{lk} : l,k \ne 0 \}$ with coefficients in $\C[\tilde{X}]$.
\begin{lem}
\label{lem:scale}
Let $\tilde{x} \in \tilde{X}$ and $U\subset X$ be a neighborhood of $\pi(\tilde{x})$
with coordinates $\{ \xi_1, \ldots, \xi_n\}$ satisfying
$U\cap C_{\tilde{x}} = \{ \xi_2 = \ldots = \xi_n = 0 \}$ and $\xi_j(\pi(\tilde{x})) = 0$ $\forall j$.

If $Z_{\tilde{x}} \cap X = \emptyset$, then there exists  a neighborhood 
$\tilde{U} \subset \tilde{X}$ of $\tilde{x}$ such that
the restriction of $\xi_1$ to $\tilde{p}(\tilde{y}) \cap U$
is an unramified covering of degree $n_{\tilde{x}}$ for all $y \in \tilde{U}$.
\end{lem}
\begin{proof}
The local coordinates define a scale, that is adapted to $\tilde{p}(\tilde{x})$, since
$U\subset X$ and $X\cap \tilde{p}(\tilde{x}) = n_{\tilde{x}} C_{\tilde{x}}$ (\cite{B}, \cite{M}).
The resulting covering is unramified, because distinct $\C$-orbits do not intersect.
\end{proof}
\begin{cor}
\label{cor:proper}
The $\C$-action on $X$ is proper if and only if
$Z_{\tilde{x}} \cap X = \emptyset$ and $n_{\tilde{x}} = 1$ for all $\tilde{x}\in\tilde{X}$.
\end{cor}
\begin{proof}
If the $\C$-action on $X$ is proper, $Z_{\tilde{x}} \cap X =\emptyset$, because the quotient
topology of the action is Hausdorff.
Furthermore, there is a local slice $S$ through $\pi(\tilde{x})$.
By construction we can assume that $S$ is a locally closed complex submanifold of
dimension $n-1$.
Choosing local coordinates $\{Ê\xi_2,\ldots \xi_n \}$ on $S$ around $\pi(\tilde{x})$ and a 
coordinate $\xi_1$ on $\C$
the biholomorphic map $\Phi: \C \times S \to \C.S$, $(t,s) \mapsto (t.s)$, gives
local coordinates on a small neighborhood $U$ of $\pi(\tilde{x})$ with the properties required
in Lemma \ref{lem:scale}.
Since $\Phi$ is a local trivialization of the action, it then follows that $n_{\tilde{x}} = 1$.

Now, let us assume $Z_{\tilde{x}} \cap X =\emptyset$ and $n_{\tilde{x}} = 1$.
Choosing coordinates $\{ \xi_1, \ldots, \xi_n\}$ as described in Lemma \ref{lem:scale} 
the locally compact complex submanifold $\pi(\tilde{U}) \cap \{ \xi_1 = 0\}$ 
is a local slice through $\pi(\tilde{x})$.
\end{proof}

\section{Examples}
\label{sec:examples}
In this section we apply the explicit formula for the map $p:\to {\mathcal C}_{d,n}$ and the
conditions derived in Corollary \ref{cor:hausdorff} and Corollary \ref{cor:proper}
to free $\C$-actions of degree two.  
These examples have been examined by different authors using other methods.
An expansion of the determinant that defines $P$ in equation (\ref{eq:P(x)})
using $f_{l k} =-f_{k l}$ yields $P = 2f_{1 0}f_{2 1}-f_{2 0}^2$ if $d=2$.

\begin{example}
% $\C^2$-action on $\C^4$, algebraic geometric quotient isomorphic to $\C^2$
We consider the free $\C$-action on $\C^3$ that is defined by the derivation
$ \delta = -x_1^2 \frac{\partial}{\partial x_2} + (1-x_1x_2) \frac{\partial}{\partial x_3}$.
It is of degree two and $U_0 = \{ x_1 \ne 0 \}$, since
$\delta^2 = x_1^3 \frac{\partial}{\partial x_3}$.
We know that this action, which is induced by a free affine $\C^2$-action on $\C^4$, 
is proper and admits a global algebraic slice \cite{Kali}, \cite{Snow}. 
Applying formula (\ref{eq:flk}) we obtain
$f_{10}=\sum_{j=0}^3((1-x_1x_2)\gamma_{3j}-\gamma_{2j}x_1^2)x_j$,
$f_{20}=\sum_{j=0}^2 x_1^3x_j\gamma_{3j}$ and  $f_{21}=-x_1^5\gamma_{32}$.
Then by equation (\ref{eq:P(x)})
\[ p(x) = -x_1^5 ( 2\gamma_{32}\sum_{j=0}^3((1-x_1x_2)\gamma_{3j}-\gamma_{2j}x_1^2)x_j
	+x_1(\sum_{j=0}^2 x_j\gamma_{3j})^2) . \]
The map $p$ extends to a regular map on $X$ with
$p|_{\{ x_1=0\}} = \gamma_{32}(\gamma_{30}x_0 +\gamma_{32}x_2)$.
For $\tilde{x} \in X\setminus U_0 $ 
the first factor corresponds the component $Z_{\tilde{x}}$, which lies at infinity, 
and the second factor represents the component $C_{\tilde{x}}$, i.e., $n_{\tilde{x}} = 1$.
\end{example}
\begin{example}
%Winkelmann
The derivation $\delta=(x_1-x_2x_5)\frac{\partial}{\partial x_3}+x_2\frac{\partial}{\partial x_4}
	+(x_1+1)\frac{\partial}{\partial x_5}$ 
defines a free $\C$-action on $\C^5$. 
Since $\delta^2=-x_2(x_1+1)\frac{\partial}{\partial x_3}$, we have $U_0 = \{ x_2(x_1+1) \ne 0 \}$
and $d=2$.
We know that this action is proper and admits Zariski-open local trivializations \cite{Wink}. 
Using the equations (\ref{eq:P(x)}) and (\ref{eq:flk}) we get
\[ p(x) = 2((x_1+1)\gamma_{35}+x_2\gamma_{34})
	(\sum_{j=0}^5 ((x_1-x_2x_5)\gamma_{3j}+x_2\gamma_{4j}+(x_1+1)\gamma_{5j})x_j )
	+x_2(x_1+1)r \]
with $r=(\sum_{j=0}^4\gamma_{3j}x_j)^2$.
The map $p:U_0 \to \overline{{\mathcal C}_{2,5}}$ extends to a regular map on 
$\{ x_2 \ne 0 \} \cup \{Êx_1+1 \ne 0 \}$ and
\begin{align*} 
p|_{\{ x_1+1=0\}}  & = -x_2\gamma_{34} 
	(\sum_{j=0}^5 ((1+x_2x_5)\gamma_{3j}-x_2\gamma_{4j})x_j ) \\
p|_{\{ x_2=0\}}  & = (x_1+1)\gamma_{35}
	(\sum_{j=0}^5 (x_1\gamma_{3j}+(x_1+1)\gamma_{5j})x_j )
\end{align*}
Hence, $Z_x \cap X = \emptyset $ and $n_x = 1$ for all 
$x \in \{ x_2 \ne 0 \} \cup \{Êx_1+1 \ne 0 \}$.
To find $\tilde{p}$ we have to construct the blow-up $\pi: \tilde{X} \to X$ 
of $X$ along $Y:=\{ x_2=x_1+1=0\}$.
There is a covering of $\tilde{X}$ by
$U_1=\{ \xi_2\ne 0,x_1+1=x_2\frac{\xi_1}{\xi_2}\}$ and
$U_2 = \{ \xi_1\ne 0,x_2=(x_1+1)\frac{\xi_2}{\xi_1}\}$.
A direct calculation shows that the rational map $p: X \to \overline{{\mathcal C}_{2,5}}$
lifts to a regular map $\tilde{p}: \tilde{X} \to \overline{{\mathcal C}_{2,5}}$ given by
\begin{align*} 
\tilde{p}|_{U_1} = & 2(\xi\gamma_{35}+\gamma_{34})
	\sum_{j=0}^5 ((x_1-x_2x_5)\gamma_{3j}+x_2\gamma_{4j}+(x_1+1)\gamma_{5j})x_j 
	+x_2 r \\
\tilde{p}|_{U_2} = & 2(\gamma_{35}+\xi\gamma_{34})
	\sum_{j=0}^5 ((x_1-x_2x_5)\gamma_{3j}+x_2\gamma_{4j}+(x_1+1)\gamma_{5j})x_j 
	+(x_1+1)r.
\end{align*}
with $r = (\sum_{j=0}^4\gamma_{3j}x_j)^2$.
Now,
\[Ê
p_1|_{U_1\cap\pi^{-1}(Y)}= (\xi\gamma_{35}+\gamma_{34})\sum_{j\ne 2,3}\gamma_{3j}x_j, 
\,
p_2|_{U_2\cap\pi^{-1}(Y)}= (\xi\gamma_{34}+\gamma_{35})\sum_{j\ne 2,3}\gamma_{3j}x_j,
\]
where the second factors correspond to $C_{\tilde{x}}$ and the first factors to $Z_{\tilde{x}}$.
\end{example}
\begin{example}
%Winkelmann2?
The derivation $\delta=x_1\frac{\partial}{\partial x_2}+x_2\frac{\partial}{\partial x_3}+
(x_2^2-2x_1x_3-1)\frac{\partial}{\partial x_4}$ defines a free $\C$-action on $\C^4$ of degree
two with $\delta^2=x_1\frac{\partial}{\partial x_3}$ and $U_0 = \{ x_1 \ne 0 \}$.
We know that this action has non-Hausdorff quotient topology \cite{Wink}.
Evaluating the map $\tilde{p}$ we can determine all inseparable orbits.
As in the examples above we find the explicit form of $p: U_0 \to {\mathcal C}_{2,4}$ using
the equations (\ref{eq:P(x)}) and (\ref{eq:flk}):
\[ p= 2 (x_1\gamma_{32}+(x_2^2-2x_1x_3-1)\gamma_{34}) 
	(\sum_{j=0}^4(x_1\gamma_{2j}+x_2\gamma_{3j}+(x_2^2-2x_1x_3-1)\gamma_{4j})x_j)
	-x_1r \]
with $r =(\sum_{j=0}^4\gamma_{3j}x_j)^2$.
This map is well defined outside of $Y:=\{ x_1=x_2^2-1=0\}$ and
\[ p|_{\{ x_1 = 0 \}} = (x_2^2-1) \gamma_{34}
	(\sum_{j=0}^4(x_2\gamma_{3j}+(x_2^2-1)\gamma_{4j})x_j) \]
which means that $n_x = 1$ and $Z_x = (x_2^2-1) \gamma_{34}$ for 
$x \in \{ x_1 = 0, x_2^2\ne 1\}$.
The blow-up $\pi: \tilde{X} \to X$ of $X$ along $Y$ is covered by the charts
$U_1^+ = \{ \xi_1\ne 0,x_2+1=x_1\frac{\xi_2}{\xi_1}\}$,
$U_2^+ = \{ \xi_2\ne 0,x_1=(x_2+1)\frac{\xi_1}{\xi_2}\}$, 
$U_1^- = \{ \xi_1\ne 0,x_2-1=x_1\frac{\xi_2}{\xi_1}\}$, and
$U_2^- = \{ \xi_2\ne 0,x_1=(x_2-1)\frac{\xi_1}{\xi_2}\}$.
It can be checked by an explicit calculation that the rational map 
$p:X \to \overline{{\mathcal C}_{2,4}}$ lifts to regular map 
$\tilde{p}: \tilde{X} \to \overline{{\mathcal C}_{2,4}}$ and
\begin{align*}
\tilde{p}|_{U_1^+\cap\pi^{-1}(Y)} = & (\gamma_{30}-\gamma_{32}+x_4\gamma_{34})
	(\gamma_{30}+\gamma_{32}+(x_4-4x_3-4\xi)\gamma_{34}), \\
\tilde{p}|_{U_2^+\cap\pi^{-1}(Y)} = & (\gamma_{30}-\gamma_{32}+x_4\gamma_{34})
	(\xi\gamma_{30}+\xi\gamma_{32}+(\xi x_4-4\xi x_3-4)\gamma_{34}), \\
\tilde{p}|_{U_1^-\cap\pi^{-1}(Y)} = & (\gamma_{30}+\gamma_{32}+x_4\gamma_{34})
	(\gamma_{30}-\gamma_{32}+(x_4+4x_3-4\xi)\gamma_{34}), \\
\tilde{p}|_{U_2^-\cap\pi^{-1}(Y)} = & (\gamma_{30}+\gamma_{32}+x_4\gamma_{34})
	(\xi\gamma_{30}-\xi\gamma_{32}+(\xi x_4+4\xi x_3-4)\gamma_{34}). \\
\end{align*}
If $\tilde{x} \in \pi^{-1}(Y)$ and $\xi(\tilde{x}) \ne 0$, then $\tilde{p}(\tilde{x})$ contains
two distinct $\C$-orbits.
In particular, orbits in $ \{ x_1=x_2-1=0 \} $ can not be separated from orbits in
$\{ x_1=x_2+1=0\}$.
\end{example}

\end{document}